\documentclass[10pt]{amsart}
\usepackage{amsmath, amsthm}
\usepackage{graphicx}
\usepackage{graphics}
\usepackage{float}

\title[Structured populations: immigration and (bi)stability]{Size-structured populations: immigration, (bi)stability and the net growth rate}
\author{J\'{o}zsef Z. Farkas}
\address{Department of Computing Science and Mathematics, University of Stirling, Stirling, FK9 4LA, UK}
\email{jzf@maths.stir.ac.uk}

\subjclass{92D25, 47D06, 35B35}
\keywords{Structured population dynamics; Population inflow; Bistability; Quasicontraction semigroups; Positivity; Spectral methods; Net growth rate}

\newtheoremstyle{theorem}
{10pt} 
{10pt} 
{\sl} 
{\parindent} 
{\bf} 
{. } 
{ } 
{} 
\theoremstyle{theorem}
\newtheorem{theorem}{Theorem}

\newtheoremstyle{defi}
{10pt} 
{10pt} 
{\rm} 
{\parindent} 
{\bf} 
{. } 
{ } 
{} 
\theoremstyle{defi}

\begin{document}

\maketitle

\begin{abstract}
We consider a class of physiologically structured population models, a first order nonlinear partial differential equation equipped with a nonlocal boundary condition, with a constant external inflow of individuals. We prove that the linearised system is governed by a quasicontraction semigroup. We also establish that linear stability of equilibrium solutions is governed by a generalized net reproduction function. In a special case of the model ingredients we discuss the nonlinear dynamics of the system when the spectral bound of the linearised operator equals zero, i.e. when linearisation does not decide stability. This allows us to demonstrate, through a concrete example, how immigration might be beneficial to the population. In particular, we show that from a nonlinearly unstable positive equilibrium a linearly stable and unstable pair of equilibria bifurcates. In fact, the linearised system exhibits bistability, for a certain range of values of the external inflow, induced potentially by All\'{e}e-effect.
\end{abstract}

\section{Introduction}
Interest in the modelling and understanding of the dynamics of biological populations is old. Classical, ordinary differential equation models assume homogeneity of individuals within population classes, and therefore involve equations for total population sizes. However, individuals in biological populations differ in their physiological characteristics. Therefore vital rates, such as those of birth, death and development, vary amongst individuals. Structured population modelling tries to account for these differences. It has played a central role in mathematical population dynamics in the past decades, being applied to amongst other systems, fish, water flea, and cell populations. The mathematical modelling of 
tumour growth and cancer treatment, human demography and the epidemiology of infectious diseases may also involve structured populations.
The first step when modelling a population may be to consider a variable that allows one to divide the population into internally
homogeneous sub-populations in order to describe the dynamics as the interaction of these groups. 
This structuring variable may be discrete \cite{CUS} or continuous \cite{I,MD,WEB} and may be based on age, size, stage of development or life cycle.
Traditionally, (continuous) structured population models have been formulated as partial
differential equations, starting with the pioneering (first nonlinear) model of Gurtin \& MacCamy \cite{GM}.    
Following the lead of \cite{PR2,PR1,WEB} we successfully applied the spectral theory of compact/positive operators to formulate biologically interpretable conditions for the linear stability/instability of equilibria for age and size-structured population models (see \cite{FJ,FJ2,FH,FH2,FH3,FH4,FH5}).
Our method allowed us to formulate very general and elegant stability/instability conditions for scramble competition models in terms of 
certain associated reproduction functions. We shall note that previously, in \cite{CR}, Rorres formulated linear stability/instability conditions in terms of the derivative of a net reproduction function for positive equilibrium solutions of an age-structured model in the special case when only the fertility depends on the standing population size. We shall also point out that Rorres' analysis in \cite{CR} was not rigorous.

Although traditionally structured population models have been formulated as partial differential equations
for population densities, the recent unified approach of Diekmann et al., making use of the rich theory of
delay and integral equations, has been resulted in significant results. Most notably
the Principle of Linearized Stability has been proven in \cite{D1,D2} for a wide class of physiologically structured population models 
formulated as delay equations (or abstract integral equations).

In the recent papers \cite{FJ2,FH3} we investigated analytically the effects of the introduction of a positive influx
of newborns on the linearised dynamical behaviour of age and size-structured scramble competition models. In our main result we were able to relax our linear stability conditions given previously in \cite{FJ,FH} for models without inflow. Motivated by our new stability results we introduced the concept of a net growth rate and formulated a very general conjecture in \cite{FJ2}. One of our goals in the present work is to show that our stability results in \cite{FH} can be extended analogously in full generality to size-structured scramble competition models with (constant) inflow at birth using the recently introduced concept of the net growth rate. 

In particular, in the present paper we consider the following nonlinear first order partial differential equation. 
\begin{equation}
p_t(s,t)+\big(\gamma(s,P(t))p(s,t)\big)_s=-\mu(s,P(t))p(s,t),\quad t>0,\ s\in (0,m),\quad m<\infty.\label{first} 
\end{equation}
Equation \eqref{first} is equipped with the following nonlinear nonlocal boundary condition representing the influx of individuals into the population
\begin{equation}
\gamma(0,P(t))p(0,t)=C+\int_0^m\beta(s,P(t))p(s,t)\,ds,\quad t>0,\quad 0\le C<\infty,\label{second}
\end{equation} 
and is subject to the initial condition
\begin{equation}
p(s,0)=:p_0(s),\quad s\in (0,m).\label{third}
\end{equation}
Here, the function $p=p(s,t)$ denotes the density of individuals of size $s$ at time $t$, hence $P(t)=\int_0^m p(s,t)\,ds$ gives the total population quantity at time $t$. $\beta$, $\mu$ and $\gamma$ denote the fertility, mortality and development rates of individuals, respectively, and all of these vital rates depend on both the structuring variable $s$ and on the total population size $P(t)$.   
We make the following general assumptions on these vital rate functions:
\begin{equation*}
\mu,\,\beta\in C^1([0,m]\times [0,\infty)),\quad \beta,\mu\geq 0,\quad\gamma\in C^2([0,m]\times [0,\infty)),\quad\gamma >0.
\end{equation*}
From the biological point of view it is natural to assume a finite maximal size $m$ which in turn helps to avoid mathematical difficulties regarding the spectral analysis of the linearised system (see e.g. \cite{FH4}). $C$ represents an inflow of minimal size individuals, i.e. an inflow of individuals from an external source at birth. There are several biological motivations for having a population inflow of individuals from an external source. In a discrete model the inflow might represent artificial stocking of individuals at certain life cycle stages for example in an insect population. In \cite{HS} Smith treated a nonlinear stage-structured matrix model with stocking and he applied his results to the celebrated "LPA model for the flour beetle Tribolium" to discover some interesting phenomenon due to the stocking. In the continuous setting, in which we are interested, a natural example for such an inflow is the case of migratory fish populations (such as a number of salmon species) that lie eggs and then move on to a different habitat. Later on the newly hatched fish join a different local population. Another example underscoring the relevance of this model is the case of fisheries where stocking of newborn (minimal size) fish is practiced. In this scenario we can regard the inflow $C$ as a control parameter in our system.

Model \eqref{first}-\eqref{third} was previously studied in \cite{CS1} by Calsina and Salda\~{n}a. They established global well-posedness and also investigated the asymptotic behaviour of solutions. In particular, they showed that under certain conditions on the vital rates \eqref{first}-\eqref{third} is governed by a strongly continuous semigroup which admits a compact global attractor. Most notably, they obtained conditions on the model ingredients which guarantee convergence of solutions to a stationary distribution when the total population tends to a constant. In contrast, our paper focuses on the local asymptotic stability of equilibrium solutions of system \eqref{first}-\eqref{third} with particular regards to the effects of the population inflow on the dynamical behaviour of the system. First, we will prove the existence of a quasicontraction semigroup describing the evolution of solutions of the linearised system. We also discuss regularity properties of the semigroup which are crucial for our stability studies. Next we will give stability/instability criterions for positive equilibrium solutions of \eqref{first}-\eqref{third}.  Then, utilizing the stability results and a ``naive'' approach, we will investigate the nonlinear dynamics of the model in a special case of model ingredients. In particular we will treat the case when, in the abscence of inflow, the strictly dominant eigenvalue of the linearised operator is zero. We will then show, through a concrete example, that the introduction of the inflow induces bistability in the system. Our claims will be supported by numerical examples, as well.

\section{Equilibrium solutions, semigroup existence and linear stability}
In this section we investigate the existence and linear stability of equilibrium solutions of \eqref{first}-\eqref{third}. We note that \eqref{first}-\eqref{third} admits the trivial solution if and only if $C=0$. When we solve \eqref{first}-\eqref{second} for a positive stationary solution $p_*$ we arrive at the following equations:
\begin{align}
& p_*(s)=p_*(0)\pi(s,P_*),\label{statsol}\\
& \gamma(0,P_*)p_*(0)=C+\int_0^m\beta(s,P_*)p_*(s)\,ds=C+p_*(0)\int_0^m\beta(s,P_*)\pi(s,P_*)\,ds,\label{2statsol}
\end{align}
where
\begin{equation}
\pi(s,P_*)=\exp\left\{-\int_0^s\frac{\gamma_s(r,P_*)+\mu(r,P_*)}{\gamma(r,P_*)}\,dr\right\}.
\end{equation}
Thus, for given vital rates $\beta,\mu,\gamma$ and inflow $C$ the function $p_*$ is a positive stationary solution of \eqref{first}-\eqref{third}
if it is determined by 
\begin{equation}
p_*(s)=\frac{P_*\pi(s,P_*)}{\int_0^m\pi(s,P_*)\,ds},\quad s\in (0,m),\label{statsolution}
\end{equation}
with the positive total population quantity 
\begin{equation}
P_*=\int_0^mp_*(s)\,ds,
\end{equation}
which satisfies the equation
\begin{equation}\label{repcond}
\gamma(0,P_*)P_*=C\int_0^m\pi(s,P_*)\,ds+P_*\int_0^m\beta(s,P_*)\pi(s,P_*)\,ds.
\end{equation}
Then, if for fixed $C\in [0,\infty)$ we define the net growth rate $Q_C$ for $P\in(0,\infty)$ as follows:
\begin{equation}
Q_C(P)=\frac{1}{\gamma(0,P)}\Big(C\,P^{-1}\displaystyle\int_0^m\pi(s,P)\,ds+\displaystyle\int_0^m\beta(s,P)\pi(s,P)\,ds\Big),\label{netgrowthrate}
\end{equation}
it is easily shown that positive equilibrium solutions of \eqref{first}-\eqref{third} are in a one-to-one correspondence with positive solutions $P_*$ of equation 
\begin{equation}
Q_C(P)=1.\label{netgrowthrateeq}
\end{equation}
Next we note that
\begin{align}
\frac{1}{\gamma(0,P)}\displaystyle\int_0^m\beta(s,P)\pi(s,P)\,ds & =\displaystyle\int_0^m\frac{\beta(s,P)}{\gamma(s,P)}\exp\left\{-\displaystyle\int_0^s\frac{\mu(r,P)}{\gamma(r,P)}\,dr\right\}\,ds \nonumber \\
& =\int_0^m\beta(a,P)\exp\left\{-\int_0^a\mu(r,P)\,dr\right\}\,da,
\end{align}
by changing variables from size $s$ to age $a$ using the relation
\begin{equation}
\frac{ds}{da}=\gamma(s(a),.).\label{varchange}
\end{equation}
Hence the function $R$, defined by 
\begin{equation}
R(P)=\frac{1}{\gamma(0,P)}\displaystyle\int_0^m\beta(s,P)\pi(s,P)\,ds=Q_0(P),
\end{equation}
is the well-known inherent net reproduction function, i.e. the expected number of newborns to be produced by an individual in her lifetime when the population size is $P$.
We also notice that $CP^{-1}$ is the per capita inflow and by \eqref{varchange} we have
\begin{equation}
L=\frac{1}{\gamma(0,P)}\displaystyle\int_0^m\pi(s,P)\,ds=\int_0^m\exp\left\{-\int_0^a\mu(r,P)\,dr\right\}\,da,
\end{equation}
where $L$ is in turn the expected lifetime of an individual. Therefore \eqref{netgrowthrate} can be rewritten as
\begin{equation}
Q_C(P)=\displaystyle\int_0^m\Big(C\,P^{-1}+\beta(a,P)\Big)\exp\left\{-\int_0^a\mu(r,P)\,dr\right\}\,da.\label{netgrowthrate2}
\end{equation}
It is straightforward to verify that the rather natural assumption
\begin{equation}
\lim_{P\to\infty}R(P)=0,
\end{equation}
yields a sufficient condition for the existence of at least one positive solution of equation \eqref{netgrowthrateeq} and in turn for the existence of 
at least one positive equilibrium solution $p_*$ of \eqref{first}-\eqref{third} via formula \eqref{statsolution}.

Given a positive stationary solution $p_*$ of system \eqref{first}-\eqref{third}, we introduce the perturbation $u=u(a,t)$ of $p$ by making the ansatz $p=u+p_*$. Then we are using Taylor series expansions of the vital rates to arrive at the linearised problem
\begin{align}
    &u_t(s,t)+ \gamma(s,P_*)\,u_s(s,t)
    +\left(\gamma_s(s,P_*)+ \mu(s,P_*)\right)\,u(s,t)\nonumber\\
    &\quad +\left(\gamma_{sP}(s,P_*)\,p_*(s)+\mu_P(s,P_*)\,p_*(s)+\gamma_P(s,P_*)\,p_*´^\prime(s)
    \right)\,\overline{U}(t)=0,\label{lin1}\\
    &u(0,t)=\int_0^m\left(\frac{\beta(s,P_*)-\gamma_P(0,P_*)p_*(0)+\int_0^m\beta_P(r,P_*)\,p_*(r)\,dr}{\gamma(0,P_*)}\right)
    \,u(s,t)\,ds\label{lin2}
\end{align}
where we have set
\begin{equation}
    \overline{U}(t)=\int_0^m u(s,t)\,ds.
\end{equation}
Eqs.~\eqref{lin1}--\eqref{lin2} are accompanied by the initial condition
\begin{equation}\label{lin3}
    u(s,0)=u_0(s).
\end{equation}
Linear semigroup theory proved to be a usefool tool to study asymptotic behaviour of solutions of similar types of physiologically structured population models, see e.g. \cite{FJ2, FH,FH2,FH3,FH4,FH5}. The first crucial step in such an analysis is to establish the existence of the governing linear semigroup. In \cite{FH,FH2,FH3} we used the Desch-Schappacher perturbation theorem to establish the existence of the governing semigroup. More recently, in \cite{FH4,FH5} we utilised dissipativity calculations to establish the existence of a governing quasicontraction semigroup. 
Here we present another approach, which is a little bit of a mix of the above two but uses the Trotter product formula, and it is potentially more useful since further regularity properties of the semigroup, such as positivity, may be immediately established. 
 
\noindent To this end, we cast the linearised system \eqref{lin1}-\eqref{lin3} in the form of an abstract Cauchy problem on the biologically relevant  state space $\mathcal{X}=L^1(0,m)$ as follows:
\begin{equation}\label{abstr}
    {d\over dt}\, u = \left({\mathcal A} + {\mathcal B}+ {\mathcal C}\right)\,u,\quad    u(0)=u_0,
\end{equation}
where
\begin{equation}
\begin{aligned}
& {\mathcal A} u = -\frac{d}{ds}\left(\gamma(\cdot,P_*)\,u\right)\quad \text{with domain}\quad\text{Dom}({\mathcal A}) =\left\{u\in W^{1,1}(0,m)\,|\,u(0)=\Phi(u)\right\},\nonumber \\
& {\mathcal B} u = -\mu(\cdot,P_*)\,u\quad \text{on ${\mathcal X}$,} \nonumber \\
& {\mathcal C} u = -\left(\gamma_{sP}(\cdot,P_*)\,p_*+\mu_P(\cdot,P_*)\,p_*+\gamma_P(\cdot,P_*)\,p_*^\prime\right)\,\int_0^m u(s)\,ds \nonumber \\
& \quad\,=\rho_*(.)\int_0^mu(s)\,ds\quad \text{on ${\mathcal X}$,} \\
& \Phi(u)= \int_0^m\left(\frac{\beta(s,P_*)-\gamma_P(0,P_*)p_*(0)+\int_0^m\beta_P(r,P_*)\,p_*(r)\,dr}{\gamma(0,P_*)}\right)\,u(s)\,ds\quad \text{on ${\mathcal X}$}.
\end{aligned}
\end{equation}
Our objective here to establish that $\mathcal{A+B+C}$ is a generator of a positive semigroup is to apply the following result (see e.g. Sect. 5 in \cite{NAG}), which is a variant of the Trotter product formula. 
\begin{theorem}\label{Trotter}
Let $\mathcal{D}$ and $\mathcal{F}$ be generators of the strongly continuous semigroups $\left\{\mathcal{T}(t)\right\}_{t\ge 0}$ and 
$\left\{\mathcal{S}(t)\right\}_{t\ge 0}$, respectively, satisfying 
\begin{equation}\label{trot}
\left|\left|\left[\mathcal{T}\left(\frac{t}{n}\right)\mathcal{S}\left(\frac{t}{n}\right)\right]^n\right|\right|\le M\exp\{\omega t\}\quad \text{for all} \quad t\ge 0, n\in \mathbb{N},
\end{equation}
for some constants $M\ge 1,\,\omega\in\mathbb{R}$. Assume that $\mathcal{Y}=\text{Dom}(\mathcal{D})\cap\text{Dom}(\mathcal{F})$ is dense in $\mathcal{X}$ 
and for some $\lambda_0>\omega$ $(\lambda_0-\mathcal{D}-\mathcal{F})\mathcal{Y}$ is also dense in $\mathcal{X}$. Then the closure of $\mathcal{D}+\mathcal{F}$ generates a strongly continuous semigroup $\{\mathcal{U}(t)\}_{t\ge 0}$ given by the formula
\begin{equation}\label{trot2}
\mathcal{U}(t)x=\lim_{n\to\infty}\left[\mathcal{T}\left(\frac{t}{n}\right)\mathcal{S}\left(\frac{t}{n}\right)\right]^nx,\quad x\in \mathcal{X}.
\end{equation}
\end{theorem}
Our first goal is to establish a condition which guarantees that $\mathcal{A}$ generates a positive quasicontraction semigroup.
We recall (see e.g. \cite{AGG}) that the linear operator $\mathcal{L}$ is called dispersive if for every $y\in\text{Dom}(\mathcal{L})\subseteq\mathcal{Y}$ there exists a $\phi_{y}\in\mathcal{Y}^*$ with $0\le\phi_y$,  $||\phi_y||_{\mathcal{Y}^*}\le 1$ and $\langle y,\phi_y\rangle_-=||y^+||_{\mathcal{Y}}$ such that 
\begin{equation*}
\langle\mathcal{L}y,\phi_y\rangle_-\le 0, 
\end{equation*}
where $\langle\cdot,\cdot\rangle_-$ stands for the natural pairing between elements of $\mathcal{Y}$ and its dual $\mathcal{Y}^*$. 
Our aim is to apply a variant of the Lumer-Phillips Theorem, see e.g. \cite{AGG,NAG}.

\begin{theorem}\label{gen}
The operator ${\mathcal A}$ generates a strongly continuous and positive quasicontraction semigroup $\{{\mathcal T_A}(t)\}_{t\geq 0}$ of bounded linear operators on ${\mathcal X}$ if the following condition holds true:
\begin{equation}\label{poscond1}
\beta(s,P_*)-\gamma_P(0,P_*)p_*(0)+\int_0^m\beta_P(r,P_*)\,p_*(r)\,dr\ge 0,\quad \text{for}\quad s\in [0,m].
\end{equation}
\end{theorem}
\begin{proof}
We consider the operator $\mathcal{A}-\omega\mathcal{I}$ for some $\omega\in\mathbb{R}$. For a fixed $u\in\mathcal{X}$ we define $\phi_u\in\mathcal{X}^*$ by $\phi_u(s)=1$ if $u(s)\ge 0$ and $\phi_u(s)=0$ if $u(s)<0$. For $u\in\text{Dom}(\mathcal{A})$ we obtain the following estimate:
\begin{align}
& ((\mathcal{A}-\omega\mathcal{I})u,\phi_u)_-\nonumber \\
& =-\int_0^m\frac{d}{ds}\left(\gamma(s,P_*)u(s)\right)\chi_{u^+}(s)\,ds-\omega\int_0^m u(s)\chi_{u^+}(s)ds \nonumber \\
& \le -\gamma(m,P_*)u(m)\chi_{u^+}(m)+\gamma(0,P_*)u(0)\chi_{u^+}(0)-\omega ||u^+||_1 \nonumber \\
& =\,b||u^+||_1-\omega||u^+||_1 \nonumber \\
& < 0,\label{estimates}
\end{align}
for some $\omega>b$ large enough, where 
\begin{equation*}
b\le \sup_{s\in[0,m]}\left\{\beta(s,P_*)-\gamma_P(0,P_*)p_*(0)+\int_0^m\beta_P(r,P_*)\,p_*(r)\,dr\right\}<\infty.
\end{equation*}
Hence we have showed that for some $\omega \in\mathbb{R}$ large enough the operator $\mathcal{A}-\omega\mathcal{I}$ is dispersive.

To obtain the range condition, i.e. that $\text{Rg}(\lambda\mathcal{I}-\mathcal{A})=\mathcal{X}$ we obtain the unique solution of the equation
\begin{equation}
(\lambda\mathcal{I}-\mathcal{A})u=h,
\end{equation}
for $\lambda$ large enough and $h\in L^1(0,m)$ as follows:
\begin{align}
u(s)= & \exp\left\{-\int_0^s\frac{\lambda+\gamma_s(y,P_*)}{\gamma(y,P_*)}\,dy\right\} \nonumber \\
 & \times\left(u(0)+\int_0^s\exp\left\{\int_0^y\frac{\lambda+\gamma_s(r,P_*)}{\gamma(r,P_*)}\,dr\right\}\frac{f(y)}{\gamma(y,P_*)}\,dy\right).\label{rangeeq}
\end{align}
Then we have
\begin{align}
|u'(s)|\le & \left|\frac{h(s)}{\gamma(s,P_*)}\right|+\left|u(0)\right|\left|\frac{\lambda+\gamma_s(s,P_*)}{\gamma(s,P_*)}\right|\exp\left\{-\int_0^s\frac{\lambda+\gamma_s(y,P_*)}{\gamma(y,P_*)}\,dy\right\}\nonumber \\
& +\int_0^s\exp\left\{-\int_y^s\frac{\lambda+\gamma_s(r,P_*)}{\gamma(r,P_*)}\,dr\right\}\frac{|h(y)|}{\gamma(y,P_*)}\,dy.
\end{align}
That is $u\in W^{1,1}(0,m)$ and the range condition is satisfied and the proof is now completed on the grounds of Theorem 1.2 in Section C-II in \cite{AGG}. 
\end{proof}
\begin{theorem}\label{generation}
Assume that condition \eqref{poscond1} and the following condition hold true
\begin{equation}\label{poscond2}
\rho_*(s)\ge 0,\quad\text{for}\quad s\in [0,m].
\end{equation}
Then, $\mathcal{A+B+C}$ generates a positive quasicontraction semigroup.
\end{theorem}
\begin{proof}
If condition \eqref{poscond2} holds true then $\mathcal{C}$ is a positive operator. 
Since $\mathcal{B+C}$ is bounded and linear and $\mathcal{B}$ is a multiplication operator it follows that 
$\mathcal{B+C}$ is a generator of a positive semigroup if condition \eqref{poscond2} holds true. Finally, since both $\mathcal{A}$ and 
$\mathcal{B+C}$ are generators of positive semigroups, Theorem 1 shows that $\mathcal{A+B+C}$ is a 
generator of a positive quasicontraction semigroup if conditions \eqref{poscond1} and \eqref{poscond2} hold true.
\end{proof}
Next we would like to establish conditions which guarantee that the equilibrium solution $p_*$ is linearly asymptotically (un)stable.
In particular, we would like to give a complete generalization of our stability results in \cite{FH} for model \eqref{first}-\eqref{third} 
with boundary immigration. To formulate stability/instability theorems one usually deduces a characteristic function (if possible!) which roots are the eigenvalues of the linearised operator $\mathcal{A+B+C}$. Then, for instability one only needs to establish the existence of a positive root of this characteristic function. 
For stability however, one also needs to guarantee that the growth bound of the semigroup, and hence stability, is governed by a real eigenvalue.
\begin{theorem}\label{regularity}
The spectrum of $\mathcal{A+B+C}$ can contain only isolated eigenvalues of finite algebraic multiplicity.
\end{theorem}
\begin{proof}
Our aim is to show that the resolvent operator $R(\lambda,\mathcal{A+B+C})$ is compact. Since $\mathcal{B+C}$ is bounded, it is enough to show that $R(\lambda,\mathcal{A})$ is compact. We obtain the solution of the resolvent equation 
\begin{equation}\label{reseq}
(\lambda-\mathcal{A})^{-1}f=u
\end{equation}
as
\begin{align}
u(s)= & \exp\left\{-\int_0^s\frac{\lambda+\gamma_s(y,P_*)}{\gamma(y,P_*)}\,dy\right\} \nonumber \\
 & \times\left(u(0)+\int_0^s\exp\left\{\int_0^y\frac{\lambda+\gamma_s(r,P_*)}{\gamma(r,P_*)}\,dr\right\}\frac{f(y)}{\gamma(y,P_*)}\,dy\right).\label{reseq2}
\end{align}
We note that, for $\lambda$ large enough the solution $u$ of the resolvent equation \eqref{reseq} belongs to $W^{1,1}(0,m)$. That is for $\lambda$ large enough the resolvent operator $(\lambda\mathcal{I}-\mathcal{A})^{-1}$ is a bounded linear mapping from $L^1(0,m)$ into $W^{1,1}(0,m)$. Since $W^{1,1}(0,m)$ is compactly embedded in $L^1(0,m)$ by the Rellich-Kondrachov theorem (see e.g. \cite{Adams}), the claim follows. 
\end{proof}
We are now able to formulate our main stability result.
\begin{theorem}\label{linstability}
Let $p_*$ be a non-trivial stationary solution with corresponding population quantity $P_*$. If 
\begin{equation}
Q'_C(P_*)>0\label{instabcond}
\end{equation} 
then $p_*$ is linearly unstable. (Note that $Q'_C$ stands for the derivative of the function $Q_C$ with respect to $P$.) 
On the other hand, suppose that conditions \eqref{poscond1}-\eqref{poscond2} and the following condition hold true:
\begin{align}
& \int_0^m p_*(s)\int_0^s\left[\frac{\gamma_{sP}(r,P_*)+\mu_P(r,P_*)}{\gamma(r,P_*)} \right. \nonumber \\
& \quad \quad \quad\quad \quad\quad\left.-\frac{\gamma_P(r,P_*)(\gamma_s(r,P_*)+\mu(r,P_*))}{\gamma^2(r,P_*)}\right]\,dr\,ds\geq -1.\label{posstrict3}
\end{align}
Then, $p_*$ is linearly asymptotically stable if 
\begin{equation}
Q'_C(P_*)<0.
\end{equation}
\end{theorem}
\begin{proof}
It is easily shown (see e.g. \cite{FH}) that the eigenvalue equation
\begin{equation}\label{eigv}
\left(\mathcal{A+B+C}-\lambda\mathcal{I}\right)u=0,\quad u(0)=\Phi(u),
\end{equation}
admits a solution $\lambda\in\mathbb{C}$ with a corresponding non-zero eigenvector $u$, if and only if $\lambda$ satisfies the following characteristic equation:
\begin{align}
1=K(\lambda)= & \int_0^mf(\lambda,s)b(s)\,ds\left(1-\int_0^mf(\lambda,s)\int_0^sf^{-1}(\lambda,y)\frac{\rho_*(y)}{\gamma(y,P_*)}\,dy\,ds\right)\nonumber \\
& +\int_0^mf(\lambda,s)\int_0^sf^{-1}(\lambda,y)\frac{\rho_*(y)}{\gamma(y,P_*)}\,dy\,ds\nonumber \\
& +\int_0^mf(\lambda,s)\,ds\int_0^mf(\lambda,s)b(s)\int_0^sf^{-1}(\lambda,y)\frac{\rho_*(y)}{\gamma(y,P_*)}\,dy\,ds,\label{characteristic}
\end{align}
where we defined the following functions:
\begin{align*}
f(\lambda,s)= & \exp\left\{-\int_0^s\frac{\lambda+\gamma_s(y,P_*)+\mu(y,P_*)}{\gamma(y,P_*)}\,dy\right\},\quad s\in [0,m],\,\lambda\in\mathbb{C}, \\
b(s)= & \frac{\beta(s,P_*)-\gamma_P(0,P_*)p_*(0)+\int_0^m\beta_P(r,P_*)\,p_*(r)\,dr}{\gamma(0,P_*)},\quad s\in [0,m].
\end{align*}
It follows from the growth behaviour of the function $f$ introduced above that the characteristic function satisfies 
\begin{equation*}
\lim_{\lambda\to +\infty}K(\lambda)=0.
\end{equation*}
Straightforward claculations yield:
\begin{align}
K(0)= & -\int_0^m\pi(s,P_*)b(s)\,ds\int_0^m\pi(s,P_*)\int_0^s\frac{\rho_*(y)}{\pi(y,P_*)\gamma(y,P_*)}\,dy\,ds \nonumber \\
 & +\int_0^m\pi(s,P_*)b(s)\,ds+\int_0^m\pi(s,P_*)\int_0^s\frac{\rho_*(y)}{\pi(y,P_*)\gamma(y,P_*)}\,dy\,ds \nonumber \\
 & +\int_0^m\pi(s,P_*)\,ds\int_0^m\pi(s,P_*)b(s)\int_0^s\frac{\rho_*(y)}{\pi(y,P_*)\gamma(y,P_*)}\,dy\,ds \nonumber \\
 = & \int_0^m\pi(s,P_*)b(s)\,ds+P_*\int_0^m\pi_P(s,P_*)b(s)\,ds+p_*(0)\int_0^m\pi_P(s,P_*)\,ds \nonumber \\
 & -p_*(0)\int_0^m\pi_P(s,P_*)\,ds\int_0^m\pi(s,P_*)b(s)\,ds,\label{calc1}
\end{align}
where $\pi_P$ stand for the derivative of the function $\pi$ with respect to $P$ and we used the relation
\begin{equation*}
P_*=p_*(0)\int_0^m\pi(s,P_*)\,ds.
\end{equation*}
On the other hand, we have 
\begin{align}
Q'_C(P_*)=R'(P_*) & +\frac{\frac{d}{dP}\left(\int_0^m\pi(s,P)\,ds\right)|_{P=P_*}\times C\gamma(0,P_*)P_*}{(\gamma(0,P_*)P_*)^2} \nonumber \\
& -\frac{C(\gamma_P(0,P_*)P_*+\gamma(0,P_*))\int_0^\pi(s,P_*)\,ds}{(\gamma(0,P_*)P_*)^2},\label{calc2}
\end{align}
and
\begin{align}
&  R'(P_*)\gamma(0,P_*)  \nonumber \\
 & = \int_0^m\left(\beta_P(s,P_*)\pi(s,P_*)+\beta(s,P_*)\pi_P(s,P_*)\right)\,ds-R(P_*)\gamma_P(0,P_*).\label{calc3}
\end{align}
Further algebraic manipulations of equations \eqref{calc2}-\eqref{calc3} yield:
\begin{align}
P_*Q'_C(P_*)= & \frac{P_*}{\gamma(0,P_*)}\left(\int_0^m\beta_P(s,P_*)\pi(s,P_*)\,ds+\int_0^m\beta(s,P_*)\pi_P(s,P_*)\,ds\right) \nonumber \\
& -R(P_*)P_*\frac{\gamma_P(0,P_*)}{\gamma(0,P_*)}+\frac{C}{\gamma(0,P_*)}\int_0^m\pi_P(s,P_*)\,ds \nonumber \\
& -C\int_0^m\pi(s,P_*)\,ds\left(\frac{\gamma_P(0,P_*)}{\gamma^2(0,P_*)}+\frac{1}{P_*\gamma(0,P_*)}\right). \label{calc4}
\end{align}
Comparing equations \eqref{calc1} and \eqref{calc4} we arrive at
\begin{equation}
K(0)=P_*Q'_C(P_*)+1.\label{calc5}
\end{equation}
Therefore, if condition \eqref{instabcond} holds true then from the Intermediate Value Theorem it follows that there exists a positive eigenvalue 
and the stationary solution $p_*$ is linearly unstable.

Next we observe that condition \eqref{posstrict3} guarantees that the characteristic function $K(\lambda)$ defined in \eqref{characteristic} 
is monotone decreasing for $\lambda\ge 0$. Therefore, equation \eqref{calc5} shows that $\mathcal{A+B+C}$ does not have a non-negative real eigenvalue 
if $Q'_C(P_*)<0$ holds. 
Next we recall that for the strongly continuous semigroup $\left\{\mathcal{T}(t)\right\}_{t\ge 0}$ with generator $\mathcal{A+B+C}$ we have 
\begin{equation*}
\omega_0=max\left\{s(\mathcal{A+B+C}),\omega_{ess}(\mathcal{A+B+C})\right\},
\end{equation*}
see e.g. \cite{NAG}, where $\omega_{ess}$ stands for the essential growth bound of the semigroup. However, for a positive semigroup the spectral bound belongs to the spectrum unless the spectrum is empty. Theorem 4 implies that the essential spectrum of the generator $\mathcal{A+B+C}$ is empty, 
therefore the growth bound is determined by a leading eigenvalue unless the spectrum is empty. 
In either case, we obtain for the growth bound $\omega_0$ of the linearised semigroup 
\begin{equation*}
\omega_0=s(\mathcal{A+B+C})<0,
\end{equation*}
and the proof of the theorem is completed.
\end{proof}

\section{Nonlinear dynamics and bistability}

In \cite{FH} our aim was to formulate linear stability/instability conditions for the model \eqref{first}-\eqref{third} in the abscence of inflow. 
Therefore, the case $R'(P_*)=Q'_0(P_*)=0$ was left completely untreated. In this case, zero is an eigenvalue of the linearised operator, see equation \eqref{calc5}. If the positivity conditions \eqref{poscond1}-\eqref{poscond2} also hold true, then zero is in fact a strictly dominant eigenvalue, hence linearisation does not decide stability. The goal of the present section is to consider this interesting and rather difficult case. To our knowledge there are a very few results in the literature (see \cite{CHQ}) which treat the case when the spectral bound of the linearised operator equals zero. We point out that our discussion will be less rigorous and some of the details need to be elaborated. For basic concepts and results used throughout this chapter we refer the reader to \cite{AGG,Car,NAG}. 

For the rest of the section we assume that conditions \eqref{poscond1},\eqref{poscond2} and \eqref{posstrict3} hold true. The semigroup $\{\mathcal{T}(t)\}_{t\ge 0}$ generated by the operator $\mathcal{A+B+C}$ is eventually compact. The spectrum of the linearised operator $\mathcal{A+B+C}$, consists of isolated eigenvalues of finite (algebraic) multiplicity. Therefore we may split the spectrum $\sigma(\mathcal{A+B+C})$ into the union of two disjoint closed subsets as follows:
\begin{equation}
\sigma(\mathcal{A+B+C})=\{\,0\}\cup\Big(\sigma(\mathcal{A+B+C})\setminus\{\,0\}\Big).\label{spectralcomp}
\end{equation}
Moreover, $\mathcal{X}$ splits into the direct sum of two  $\{\mathcal{T}(t)\}_{t\ge 0}$-invariant closed subspaces: $\mathcal{Y}$ and $\mathcal{Z}$. In fact, there exists a spectral decomposition of $\mathcal{X}$ with respect to \eqref{spectralcomp}. This is because the set $\{\,0\}$ in \eqref{spectralcomp} is trivially bounded, see e.g. \cite{AGG}. Note that, in general, the set $\sigma(\mathcal{A+B+C})\setminus\{\,0\}$ is non-empty, hence the spectral decomposition is not trivial, but the formulation of precise conditions is left for future work. Moreover, we may determine the subspace $\mathcal{Y}$ explicitly. Straightforward calculations show that the eigenvalue equation 
\begin{equation}
(\mathcal{A+B+C})U=0
\end{equation}
subject to 
\begin{equation}
U(0)=\int_0^m\left(\frac{\beta(s,P_*)-\gamma_P(0,P_*)p_*(0)+\int_0^m\beta_P(r,P_*)\,p_*(r)\,dr}{\gamma(0,P_*)}\right)\,U(s)\,ds,
\end{equation}
has a one-dimensional solution space spanned by the function
\begin{equation}
F(s)\left(\frac{1+\int_0^mF(s)\int_0^s\frac{G(r)}{F(r)}\,dr\,ds}{\int_0^mF(s)\,ds}-\int_0^s\frac{G(r)}{F(r)}\,dr\right)\in L^1(0,m),
\end{equation}
where
\begin{equation}
F(s)=\exp\left\{-\int_0^s\frac{\gamma_s(r,P_*)+\mu(r,P_*)}{\gamma(r,P_*)}\,dr\right\},\quad s\in[0,m],
\end{equation}
and
\begin{equation}
G(r)=\gamma_{sP}(r,P_*)p_*(r)+\mu_P(r,P_*)p_*(r)+\gamma_P(r,P_*)p'_*(r),\quad r\in [0,m].
\end{equation}
Hence the one-dimensional subspace $\mathcal{Y}$ of $\mathcal{X}$ given by
\begin{equation}
\mathcal{Y}=\left\{F(.)\left(\frac{1+\int_0^mF(s)\int_0^s\frac{G(r)}{F(r)}\,dr\,ds}{\int_0^mF(s)\,ds}-\int_0^.\frac{G(r)}{F(r)}\,dr\right)\right\},
\end{equation}
which is $(\mathcal{A+B+C})$-invariant, is the centre eigenspace. We note that the set $\sigma(\mathcal{A+B+C})\setminus \{0\}$ is contained in the half plane $H_\delta=\{z\in\mathbf{C}\,:\, \text{Re}(z)<\delta\}$ for some $\delta<0$. This, and the existence of the spectral decomposition implies that the semigroup $\{\mathcal{S}(t)\}_{t\ge 0}$, which is the restriction of $\{\mathcal{T}(t)\}_{t\ge 0}$ to $\mathcal{Z}$ is uniformly exponentially stable (see e.g. \cite{AGG}). There exists a centre manifold (see e.g. \cite{Car}), but it cannot be determined explicitly since the eigenvalues of $\mathcal{A+B+C}$ (and in turn the corresponding eigenspaces) are not explicitly available (except the zero eigenvalue). 

In the rest of the section we restrict our attention to the special case of vital rates:
\begin{equation}
\mu=\mu(s),\quad \gamma=\gamma(s),\quad \beta=\beta(s,P),\quad s\in [0,m],\quad P\in [0,\infty),
\end{equation}
and we further assume that 
\begin{align}
 R'(P_*) & =\int_0^m\beta_P(s,P_*)\exp\left\{-\int_0^s\frac{\gamma_s(r)+\mu(r)}{\gamma(r)}\,dr\right\}\,ds=0,\nonumber\\
 R''(P_*) & =\int_0^m\beta_{PP}(s,P_*)\exp\left\{-\int_0^s\frac{\gamma_s(r)+\mu(r)}{\gamma(r)}\,dr\right\}\,ds\ne 0.
\end{align} 
In this scenario, equation \eqref{first} is linear, and the only nonlinearity arises in the boundary condition \eqref{second}. Since linearisation does not decide stability we take into account the quadratic term when using Taylor series expansion for the fertility function $\beta$. We obtain the nonlinear problem:
\begin{align}
& \frac{d}{dt}U=(\mathcal{A+B})U,\nonumber \\
& \text{Dom}(\mathcal{A})= \nonumber \\
& \left\{U\in W^{1,1}(0,m)|U(0)=\int_0^m\frac{\beta(s,P_*)}{\gamma(0)}U(s)\,ds+\overline{U}^2\int_0^m\frac{\beta_{PP}(s,P_*)}{\gamma(0)}p_*(s)\,ds\right\},\nonumber \\
& \quad \overline{U}=\int_0^mU(s)\,ds.\label{nonlin}
\end{align}
When we look to obtain solutions of the equation
\begin{align}
& (\mathcal{A+B})U=\varepsilon,\quad \varepsilon\in\mathbf{R},\label{nonlin2}
\end{align}
we arrive at
\begin{equation}
0=\varepsilon\int_0^m\beta(s,P_*)F(s)\int_0^s\frac{1}{F(r)}\,dr\,ds+\overline{U}^2\int_0^m\beta_{PP}(s,P_*)p_*(s)\,ds,\label{cond}
\end{equation}
an equation for $\overline{U}$. Existence of a (real) solution of \eqref{cond} is shown to be neccessary and sufficient for the existence of a non-trivial solution of \eqref{nonlin2}. This shows that \eqref{nonlin2} admits solutions only for $\varepsilon>0$ or for $\varepsilon<0$
depending on the sign of $R''(P_*)$. This means that the flow generated by the linear operator $\mathcal{A+B}$ is monotone in a neighbourhood of $U\equiv 0$, which in turn implies that the equilibrium $p_*$ is unstable. 

This claim will be confirmed numerically in case of an example. When we make the following choice of vital rates: 
\begin{equation}
\mu\equiv 1,\quad \gamma\equiv 1,\quad \beta=\frac{P^2e^{-P}se^{-s}+0.5P^2e^{-P}}{3e^{-2}-2e^{-8}-13e^{-14}}\quad s\in [0,6], \label{example}
\end{equation}
system \eqref{first}-\eqref{third} admits for $C=0$ a unique positive equilibrium solution 
\begin{equation}
p_*(s)=2\frac{e^{-s}}{1-e^{-6}},\quad s\in [0,6],\label{posstatex}
\end{equation} 
with total population size $P_*=2$. It is straightforward to confirm that conditions \eqref{poscond1}, \eqref{poscond2} and \eqref{posstrict3} hold true 
and $Q'_0(P_*)=R'(P_*)=0$, $Q''_0(P_*)=R''(P_*)<0$ since $\beta_{PP}(s,P_*)<0$ holds for  $s\in [0,6]$. For $C=0$ we have $Q_0(0)=0<1$ which implies that the trivial equilibrium is locally asymptotically stable (see e.g. \cite{FH}). The positive equilibrium $p_*$ given by \eqref{posstatex} is nonlinearly unstable. 

When we introduce the inflow, at $C=0$ a pair of equilibria, denoted by $p_*^1$ and $p_*^2$, bifurcates from the positive equilibrium $p_*$. We may assume that 
\begin{equation}
P_*^1=\int_0^mp_*^1(s)\,ds<P_*^2=\int_0^mp_*^2(s)\,ds.
\end{equation} 
Then, it is easily shown that $Q'_C(P^2_*)<0$ holds, and conditions \eqref{poscond1}, \eqref{poscond2} and \eqref{posstrict3} are satisfied. Hence the equilibrium $p_*^2$ is linearly stable while $p_*^1$, for which  $Q'_C(P^1_*)>0$ holds, is linearly unstable by Theorem \ref{linstability}. The third positive equilibrium, which is ``connected'' to the trivial one, satisfies also the stability conditions of Theorem \ref{linstability}, hence it remains locally asymptotically stable. The system is bistable for $C\in (0,C_*)$ for some $C_*>0$ . This bistability is induced potentially by an All\'{e}e-type effect. When $C$ reaches the critical value $C_*$ the system undergoes another bifurcation. This is illustrated in Figure 1.
\begin{figure}[H]
\begin{center}
\includegraphics[width=8cm,height=5cm]{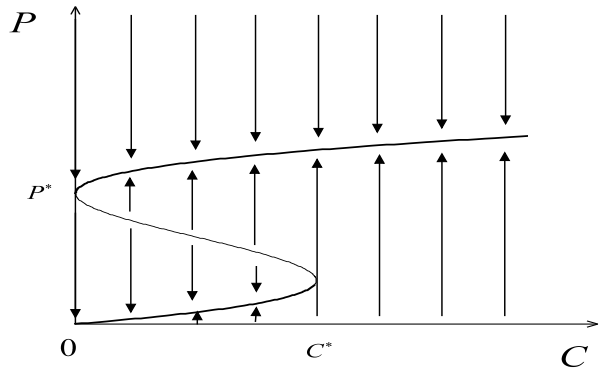}
\caption{{\it Schematic bifurcation diagram of equilibrium population sizes}}
\end{center}
\end{figure}
In Figure 2. we present the results of some numerical simulations. We used the hier-community program which solves equations \eqref{first}-\eqref{third} using a finite difference scheme approximation, see \cite{ADH}. It is developed by A.~S. Ackleh and his collaborators and it is available at Ackleh's web page.
\begin{figure}[H]
\begin{center}
\includegraphics[width=10cm,height=10cm]{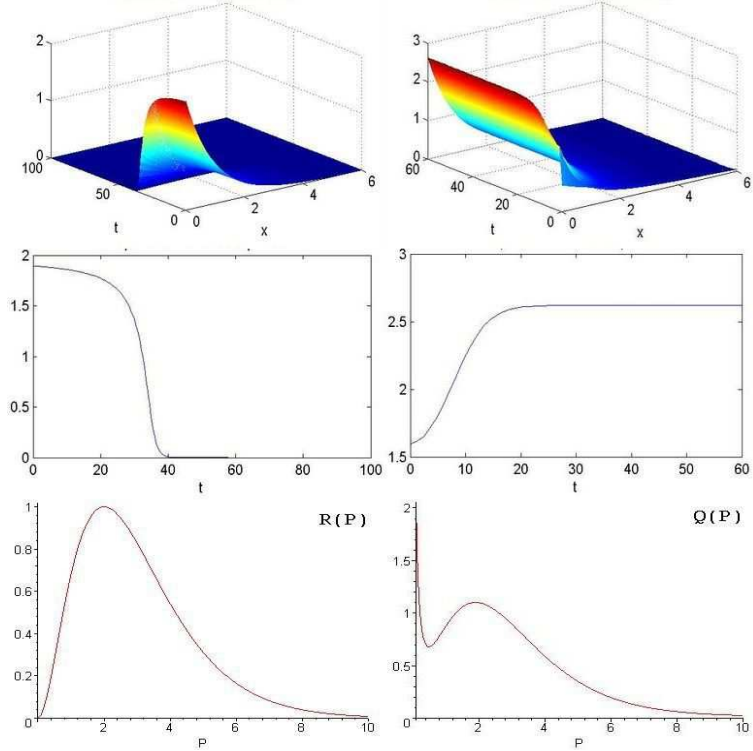}
\end{center}
\caption{{\it In the first column we plotted the evolution of the solution and its integral starting with the initial condition $p_0(x)=\frac{1.9 e^{-x}}{1-e^{-6}}$,  and the graph of the net reproduction rate $R$, respectively. (Here $x$ replaces the structuring variable $s$.) In the second column we plotted again the evolution of the solution and its integral starting with the initial population distribution $p_0(x)=\frac{0.7 e^{-0.4x}}{1-e^{-6}}$,  and the graph of net growth rate $Q$, in the case of $C=0.2$.}}
\end{figure}

\section{Concluding remarks}

In the present manuscript we employed semigroup and spectral methods to investigate the linearised dynamical behaviour of a size-structured 
population model. The model we treated here is an extension of the one analysed previously in \cite{FH}. In particular our goal was to study the effects 
of migration of individuals from an external source. This question is biologically relevant. Migration was modelled through a constant inflow of minimal size individuals. The assumption of migration at only one state might seem unsatisfactory from the biological point of view, nevertheless we managed to give a complete generalization of our stability/instability results established previously in \cite{FH} for a model without migration. In fact, in our main result we obtained biologically relevant stability/instability criteria by defining an appropriate net reproduction function of the standing population. In my opinion this result is quite interesting and it shall have an impact when studying other physiologically structured scramble competition models.

In the second half of the present work, we were mainly concerned with the weakly nonlinear dynamics of the model. 
The bifurcation diagram in Figure 1 illustrates the dynamic behaviour of the system in case of migration. 
We also gave an example when the introduction of the inflow stabilizes a weakly unstable steady state.  
Although in this part we treated a quite special case of model ingredients the analysis may well be applicable to more general scenarios. 

\bigskip

\noindent {\bf Acknowledgment}

The author was supported by the EPSRC grant EP/F025599/1.

\end{document}